\documentclass[a4paper,12pt]{article}

\usepackage{amsmath}
\usepackage{amsthm}
\usepackage{amssymb}  
\usepackage{latexsym} 
\usepackage[all]{xy}
\usepackage{comment}
\usepackage{rotating}
\usepackage{slashbox}

\newcommand{\EE}{{\mathcal E}}
\newcommand{\CC}{{\mathcal C}}
\newcommand{\Cdes}{\widetilde{\mathcal C}}
\newcommand{\mm}{\mathfrak m}
\newcommand{\length}{\operatorname{length}}

\newcommand{\Spec}{\operatorname{Spec}}
\newcommand{\level}{\operatorname{level}}

\newcommand{\GL}{\operatorname{GL}}

\newcommand{\isom}{ \cong }

\newcommand{\OK}{{\mathcal{O}_K}}

\newcommand{\Div}{\operatorname{Div}}

\newcommand{\Frac}{\operatorname{Frac}}
\newcommand{\Aff}{{\mathbb A}}
\newcommand{\PP}{{\mathbb P}}

\newcommand{\supp}{\operatorname{Supp}}
\newcommand{\Proj}{\operatorname{Proj}}

\newfont{\wncyr}{wncyr10 at 12pt}
\newfont{\wncyrten}{wncyr10 at 10pt}

\newenvironment{Proof}{\par\noindent{\sc Proof:}}%
                      {\hspace*{\fill}\nobreak$\Box$\par\medskip}
\newenvironment{ProofOf}[1]{\par\noindent{\sc Proof of #1:}}%
                       {\hspace*{\fill}\nobreak$\Box$\par\medskip}

\newenvironment{myitemize}
{\begin{itemize}
\setlength{\itemsep}{1pt}
\setlength{\parskip}{0pt}
\setlength{\parsep}{0pt}}
{\end{itemize}}

\newtheorem{Proposition}{Proposition}[section]
\newtheorem{Theorem}[Proposition]{Theorem}
\newtheorem{Lemma}[Proposition]{Lemma}
\newtheorem{Corollary}[Proposition]{Corollary}

\theoremstyle{definition}
\newtheorem{Definition}[Proposition]{Definition}
\newtheorem{Remark}[Proposition]{Remark}

\newcounter{nootje}
\begin{document}

\title{Minimal Genus One Curves}
\author{Mohammad Sadek}
\date{}
\maketitle
\begin{abstract}{\footnotesize
In this paper we consider genus one equations of degree $n$, namely a (generalised) binary quartic when $n=2$, a ternary cubic when $n=3$, and a pair of quaternary quadrics when $n=4$. A new definition for the minimality of genus one equations of degree $n$ is introduced. The advantage of this definition is that it does not depend on invariant theory of genus one curves. We prove that this definition coincides with the classical definition of minimality when $n\le4$. As an application, we give a new proof for the existence of global minimal genus one equations over number fields of class number 1.}
\end{abstract}

\section{Introduction}
\label{sec:intro}

We work throughout over a Henselian discrete valuation field $K$ with ring of integers $\OK$. We fix a uniformiser $t$, a normalisation $\nu$, and write $k=\OK/t\OK$. Set $S=\Spec\OK$. Let $C$ be a smooth genus one curve over $K$. $D$ will denote a $K$-rational divisor on $C$ of degree $n$. If $n=1$, then $C(K)\ne\emptyset$. If $n\ge 2$, then the divisor class $[D]$ defines a morphism $C\to\PP_K^{n-1}$. Let $R$ be a Dedekind domain. A {\em genus one equation of degree $n$} describes the pair $(C,D)$ and is defined as follows:

 If $n=1,$ a Weierstrass equation \begin{equation} y^2 + a_1 x y + a_3 y = x^3 + a_2 x^2 + a_4 x + a_6,\;a_i\in K.\end{equation}
Two genus one equations of degree $1$ with coefficients in $R$ are {\em $R$-equivalent} if they are related by the substitutions:
$x' = u^2x + r$ and $y' = u^3y + su^2x + t,$ where $r,s,t\in R,\;u\in R^*.$ We set $\det([u;r,s,t])=u^{-1}.$

 If $n=2,$ a (generalised) binary quartic \begin{equation}y^2 + (\alpha_0 x^2 + \alpha_1 xz + \alpha_2z^2) y = a x^4 + b x^3z + c x^2z^2 + d xz^3 + ez^4.\end{equation}Two genus one equations of degree $2$ with coefficients in $R$ are {\em $R$-equivalent} if they are related by the substitutions: $x'=m_{11}x+m_{21}z,\;z'=m_{12}x+m_{22}z$ and $y'=\mu^{-1}y+r_0x^2+r_1xz+r_2z^2,$ where $\mu\in R^*,\; r_i\in R$ and $M=(m_{ij})\in \GL_2(R)$. We set $\det([\mu,(r_i),M])=\mu\det M.$

 If $n=3,$ a ternary cubic
 \begin{equation}F(x,y,z)=a x^3 + b y^3 + c z^3 + a_2x^2 y + a_3x^2 z + b_1y^2x+ b_3 y^2z + c_1z^2x + c_2 z^2y + m x y z=0.\end{equation} Two genus one equations of degree $3$ with coefficients in $R$ are {\em $R$-equivalent} if they are related by multiplying by $\mu\in R^*$, then substituting
 $x'=m_{11}x+m_{21}y+m_{31}z,\;y'=m_{12}x+m_{22}y+m_{32}z$ and $z'=m_{13}x+m_{23}y+m_{33}z$, where $M=(m_{ij})\in \GL_3(R).$
 Set $\det([\mu,M])=\mu\det M.$

 If $n\ge4$, then $C$ is described by $l=n(n-3)/2$ quadratic forms in $n$ variables. Two genus one equations of degree $n$ with coefficients in $R$ are {\em $R$-equivalent} if they are related by the substitutions: $F'_i=m_{i1}F_1+m_{i2}F_2+\ldots + m_{il} F_l,\; i=1,\ldots,l,$ for $M=(m_{ij})\in\GL_l(R),$
and then $x'_j=\sum_{i=1}^nn_{ij}x_i$ for $N=(n_{ij})\in\GL_n(R).$
When $n=4$, set $\det([M,N])=\det M\det N.$

For $n\le4$, we associate invariants $c_{4,\phi}$, $c_{6,\phi}$ and {\em discriminant} $\Delta_{\phi}$ to a genus one equation $\phi$ of degree $n$
such that $\Delta_{\phi}=(c_{4,\phi}^3-c_{6,\phi}^2)/1728.$ Moreover, $\phi$ defines a smooth curve of genus one if and only if $\Delta_{\phi}\ne 0.$ The invariants $c_{4,\phi},c_{6,\phi}$ and $\Delta_{\phi}$ are of weights $r=4,6$ and $12$ respectively. In other words, if $F\in\{c_{4,\phi},c_{6,\phi},\Delta_{\phi}\}$, then $F\circ g= (\det g)^rF$ for every $g$ defined by an $R$-equivalence. These invariants have been known since the nineteenth century, and can be found in \cite{AKMMMP}.
We scale these invariants according to \cite{FiInvariants}.

\begin{Definition}\label{def:integral-minimal}
 A genus one equation $\phi$ of degree $n$, $n\le4$, with $\Delta_{\phi}\ne 0$ is
 \begin{myitemize}
\item[(a)] {\em integral} if the defining polynomials have coefficients in $\OK.$
\item[(b)] {\em minimal} if it is integral and $\nu(\Delta_{\phi})$ is
minimal among all the valuations of the discriminants of the integral genus one equations $K$-equivalent to $\phi.$
\end{myitemize}
\end{Definition}

Producing integral genus one equations of degree $n$ with small coefficients has been a target for investigations. In order to obtain such genus one equations, we need to {\em reduce} and {\em minimize} them. By reducing genus one equations, we mean reducing the size of the coefficients by a unimodular linear change of coordinates, which does not change the invariants. To minimise genus one equations, we need to make the associated invariants smaller. Reduced and minimal genus one equations of degree $2$ appear as an essential part of the 2-descent algorithm described by Birch and Swinnerton-Dyer in \cite{swinnerton}. More recent treatment for the minimisation of genus one equations of degree 2, 3 and 4 can be found in \cite{CremonaStollMini}, \cite{Fiminimise} and \cite{WomackThesis} respectively. An algorithmic approach for minimising genus one equations of degree $n$, $n\le4,$ can be found in \cite{FiStCr}. This paper will be dedicated to the minimisation question.

The main obstacle which holds back the existence of a neat solution for the minimsation question when $n$ is large is the difficulty of describing the rings of invariants associated to these genus one equations. In order to overcome this difficulty, we give an alternative definition for the minimality of genus one equations of degree $n$.

 If $\phi$ is an integral genus one equation of degree $n$, then it defines an $S$-scheme $\CC$. We give criteria for $\CC$ to be normal, and hence an $S$-model for its generic fiber, see Definition \ref{def:models} below. Our definition of minimality compares the minimal desingularisation of $\CC$ to the minimal proper regular model of its generic fiber. This definition does not involve invariant theory of genus one curves. We prove that our definition agrees with Definition \ref{def:integral-minimal} when $n\le 4$. In a sequel to this paper, we will use our definition of minimality to count the number of minimal genus one equations of degree $n$ up to $\OK$-equivalence in a given $K$-equivalence class, for $n\le 4$.
Finally, we give a new proof for the following theorem (\cite{FiStCr}, Theorem 4.17).

 \begin{Theorem}
 \label{thm1}
 Let $F$ be a number field of class number $1$ with ring of integers $\mathcal{O}_F$. Let $C$ be a smooth genus one curve defined over $F$ by a genus one equation $\phi$ of degree $n$ for $n\le4$. Assume that $C(F_{\nu})\ne\emptyset$ for every completion $F_{\nu}$ of $F$. Let $E$ be the Jacobian elliptic curve of $C$ with minimal discriminant $\Delta$. Then $\phi$ is $F$-equivalent to an $\mathcal{O}_F$-integral genus one equation whose discriminant is $\Delta.$
 \end{Theorem}

\section{Models of genus one curves}
\label{sec:models of genus one curves}

\begin{Definition}\label{def:models}
An {\em $S$-curve} is an integral, projective, flat, normal $S$-scheme $f:X\to S$ of dimension 2. The generic fiber of $\CC$ will be denoted by $\CC_K$ and its special fiber by $\CC_k$.
We define an {\em $S$-model} for a smooth curve $C$ over $K$ to be an $S$-curve $\CC$ such that $\CC_K$ is isomorphic to $C$.
\end{Definition}

\begin{Definition}\label{def:contraction}
Let $\CC$ be an $S$-curve. Let $(\Gamma_i)_{i\in I}$ be the family of irreducible components of $\CC_k.$ For a strict subset $J\subset I$, a \emph{contraction} of the components $\Gamma_j,\;j\in J,$ in $\CC$ consists of an $S$-morphism $u:\CC\rightarrow \CC^J $ of $S$-schemes such that\\
(a) For each $j\in J$, the image $u(\Gamma_j)$ consists of a single point $x_j\in \CC^J$, and\\
(b) $u$ defines an isomorphism $\CC-\bigcup_{j\in J}\Gamma_j\xrightarrow{\sim}\CC^J-\bigcup_{j\in J}{x_j}.$
\end{Definition}
Since $\OK$ is Henselian, the contraction $u:\CC\to\CC^J$ of the components $\Gamma_j,\;j\in J,$ exists for any strict subset $J\subset I$. Moreover, the morphism $u$ is unique up to unique isomorphism, see (\cite{Liubook}, Theorem 8.3.36 and Proposition 8.3.28).
The following theorem, Theorem $1$ of (\cite{Neron}, \S 6.7), describes the contraction morphism explicitly.

\begin{Theorem}
\label{th:contraction}
Let $\CC$ be an $S$-curve. Let $(\Gamma_i)_{i\in I}$ be the family of irreducible components of $\CC_k.$ Let $D$ be a non-trivial effective Cartier divisor on $\CC.$ Let $J$ be the set of all indices $j\in I$ such that $\supp (D)\cap \Gamma_j=\emptyset.$ Then the canonical morphism
$$u:\CC\to\CC^J:=\Proj(\bigoplus_{m=0}^{\infty}H^0(\CC,\mathcal{O}_{\CC}(m D)))$$
is the contraction of the components $\Gamma_j,j\in J,$ and $\CC^J$ is an $S$-curve.
\end{Theorem}

Let $C$ be a smooth genus one curve over $K$. Assume that $C(K)\ne\emptyset$. Let $E$ be the Jacobian elliptic curve of $C$ with minimal proper regular model $E^{min}$. Since $C\isom_K E$, the minimal proper regular model $C^{min}$ of $C$ is isomorphic to $E^{min}$. Thus we will dispense with $C^{min}$ and write $E^{min}$ instead.

The $S$-scheme $\CC$ defined by an integral genus one equation $\phi:y^2+g(x,z)y=f(x,z)$ of degree $2$
is the scheme obtained by glueing
$ \{ y^2+g(x,1)y = f(x,1) \} \subset \Aff_S^2 $ and
$ \{ v^2+g(1,u)v = f(1,u) \} \subset \Aff_S^2 $ via $x=1/u$ and $y=x^2v$. It comes with
a natural morphism $\CC \to \PP^1_S$ given on these affine pieces
by $(x,y) \mapsto (x:1)$ and $(u,v) \mapsto (1:u)$.

The $S$-scheme defined by an integral genus one equation $\phi$ of degree $n$, $n=1$ or $n\ge3$, is
the subscheme $\CC \subset \PP^{m}_S$ defined by $\phi$, where $m=2$ when $n=1$, and $m=n-1$ when $n\ge3$.

Now we give the key definition of this paper.
\begin{Definition}
\label{def:geometrically minimal}
Let $\phi$ be an integral genus one equation of degree $n,\;n\ge 1$, defining a normal $S$-scheme $\CC$. Assume moreover that $\CC_K$ is smooth and $\CC_K(K)\ne\emptyset$. Let $E^{min}$ be the minimal proper regular model of the Jacobian of $\CC_K$. Then $\phi$ is said to be {\em geometrically minimal} if the minimal desingularisation $\Cdes\to\CC$ satisfies $\Cdes\isom E^{min}$.
\end{Definition}

For the definitions of minimal proper regular models and minimal desingularisations see (\cite{Liubook}, \S 9.3).

\section{Normality}
\label{sec:normality}

In this section we prove that an $S$-scheme defined by a minimal genus one equation of degree $n,\;n\le4,$ is normal, and hence is an $S$-model for its generis fiber.
Let $\CC$ be an $S$-scheme defined by an integral genus one equation of degree $n$ where $\CC_K$ is smooth. It is known that the normality of $\CC$ implies that there are only finitely many non-regular points on $\CC,$ and all these points are closed points in the special fiber. Moreover, if $n\le 4$, then $\CC$ can be seen as a complete intersection. It follows that $\CC$ is normal if and only if $\CC$ is regular at the generic points of $\CC_k.$ The latter statement is a direct consequence of Serre's criterion for normality, see (\cite{Liubook}, Corollary 8.2.24). If $\CC_k$ is reduced, then $\CC$ is normal, see (\cite{Liubook}, Lemma 4.1.18).

\begin{Lemma}
\label{lem:regularity}
Let $(A,\mm)$ be a regular Noetherian local ring.
\begin{myitemize}
\item[(i)] Suppose that $f\in\mm\backslash\{0\}.$ Then $A/fA$ is regular if and only if $f\not\in\mm^2.$
\item[(ii)] Suppose that $I$ is a proper ideal of $A.$ Then $A/I$ is regular if and only if $I$ is generated by $r$ elements of $\mm$ which are linearly independent mod $\mm^2.$
    \end{myitemize}
\end{Lemma}
\begin{Proof}
See (\cite{Liubook}, Corollaries 4.2.12 and 4.2.15).
\end{Proof}

\begin{Lemma}
\label{lem:S' normality implies S normality}
Let $K'$ be a finite extension of $K$ with ring of integers $\mathcal{O}_{K'}$.
Let $\CC$ be an $S$-scheme. Set $S'=\Spec\mathcal{O}_{K'}$, and $\CC'=\CC\times_SS'$. If $\CC'$ is $S'$-normal, then
$\CC$ is $S$-normal.
\end{Lemma}
\begin{Proof}
That $\CC'$ is $S'$-normal means that $\mathcal{O}_{\CC',x}$ is integrally closed in $\Frac(\mathcal{O}_{\CC',x})$ for every $x\in\CC'$. Now let $x\in\CC$, and $\alpha\in\Frac(\mathcal{O}_{\CC,x})$ satisfy an integral relation for $\alpha$ over $\mathcal{O}_{\CC,x}$. The $S'$-normality of $\CC'$ implies that $\alpha\in\mathcal{O}_{\CC',x}$. Therefore, $\alpha\in\mathcal{O}_{\CC',x}\cap\Frac({\mathcal{O}_{\CC,x}})=\mathcal{O}_{\CC,x}$.
\end{Proof}

If $f(x_1,\ldots,x_n)=\sum_{i=1}^m a_ix_1^{l_{1i}}\ldots x_n^{l_{ni}}\in\OK[x_1,\ldots,x_n],$ then $\tilde{f}(x_1,\ldots,x_n)$ will denote its image in $k[x_1,\ldots,x_n].$ Moreover, $\nu(f)=\min \{\nu(a_i):1\le i\le m\}.$

Let $\CC$ be an $S$-scheme defined by an integral genus one equation $\phi$ of degree $3$ where
\begin{equation}\label{eq7}
\phi: b y^3 + f_1(x,z)y^2+f_2(x,z) y+ f_3(x,z)=0,
\end{equation}
with $f_1(x,z)=b_1x+b_3z, f_2(x,z)=a_2x^2 + m x z +c_2 z^2$ and $ f_3(x,z)=ax^3+a_3x^2z+c_1z^2x+cz^3.$ If $\CC_k$ contains an irreducible component of multiplicity-$m,\;m\ge 2,$ then we can assume without loss of generality that the defining equation of this multiplicity-$m$ component is $y=0.$ This means that $\min\{\nu(f_2),\nu(f_3)\}\ge1,\nu(f_1)=0$ when $m= 2,$ and $\min\{\nu(f_1),\nu(f_2),\nu(f_3)\}\ge1,\nu(b)=0$ when $m=3.$

\begin{Proposition}
\label{prop:normal1,2,3}
 Let $\CC$ be an $S$-scheme defined by an integral genus one equation $\phi$ of degree $n,\;n\le3.$ Assume that $\CC_K$ is smooth.
 \begin{myitemize}
\item[(i)] if $\CC_k$ consists only of multiplicity-$1$ components, then $\CC$ is normal.
\end{myitemize}
Now assume that $\CC_k$ contains an irreducible component of multiplicity greater than $1$.
\begin{myitemize}
\item[(ii)] if $n=2$ and $\phi:y^2+g(x)y=f(x)$, then $\CC$ is normal if and only if there exists $R(x)\in\OK[x]$ such that $\nu(f(x)+g(x)R(x)-R(x)^2)=1$.
\item[(iii)] if $n=3$, $\phi:F(x,y,z)=0$ is given as in equation (\ref{eq7}), and $\CC_k$ contains a multiplicity-$m$ component $\Gamma:\{y=0\}$, $m\ge2$, then $\CC$ is normal if and only if $\nu(f_3)=1$.
\end{myitemize}
\end{Proposition}
\begin{Proof}
$(i)$ Since $\CC_k$ is reduced and $\CC_K$ is smooth, the $S$-scheme $\CC$ is normal.
$(ii)$ This is Lemme 2 (c) of \cite{LiuModeles}.
$(iii)$ The maximal ideal corresponding to the generic point $\xi$ of $\Gamma$ is $\mm_{\xi}=\langle t,y\rangle$. We have $\CC$ is normal if and only if $F(x,y,z)\not\in\mm_{\xi}^2$, see Lemma \ref{lem:regularity} $(i)$. Since $\nu(f_2)\ge 1$, we have $y^3,y^2,f_2(x,z)y\in\mm_{\xi}^2.$ Therefore, $F(x,y,z)\not\in\mm_{\xi}^2$ if and only if $\nu(f_3)=1.$
\end{Proof}

Now we study the normality of an $S$-scheme $\CC$ defined by an integral genus one equation $\phi:F(x_1,x_2,x_3,x_4)=G(x_1,x_2,x_3,x_4)=0$ of degree $4$. We will assume that $\tilde{F}$ and $\tilde{G}$ are coprime. Moreover, we will assume that $\phi$ is not $\OK$-equivalent to an equation whose reduction mod $t$ is given by $x_1^2=x_2^2=0$. These two assumptions are reasonable to make because of the following lemma
which can be found in \S 2.5.1 of \cite{WomackThesis}.
\begin{Lemma}
 \label{lem:quadruple line no rational points}
  \begin{myitemize}
 \item[(i)] If $\tilde{F}$ and $\tilde{G}$ have a common factor, then $\phi$ is not minimal.
 \item[(ii)] If $\CC_k$ is defined by $x_1^2=x_2^2=0$, then either $\phi$ is not minimal, or $\CC_K(K)=\emptyset$.
     \end{myitemize}
 \end{Lemma}

We observe that (i) all the irreducible components of $\CC_k$ are defined over the residue field of a finite unramified extension of $K$, (ii) the normality of $\CC$ over the ring of integers of a finite extension of $K$ implies the normality of $\CC$ over $\OK$, see Lemma \ref{lem:S' normality implies S normality}, and (iii) the minimality of $\phi$ is stable under unramified base changes, see (\cite{FiStCr}, Theorem 3.6). Therefore, we will assume that $k$ is algebraically closed when we are finding criteria for the normality of $\CC$, see Proposition \ref{prop:normal4}, and testing the normality of $\CC$ when $\CC$ is minimal, see Theorem \ref{thm:minimal+soluble imlies normal}.

Since we will be interested in $\CC$ when $\CC_k$ contains a component of multiplicity-$m$, $m\ge 2$, we will write down all the possibilities for such a special fiber, up to $\OK$-equivalence of $\phi$. Again $k$ will be algebraically closed for the remainder of this section. For a complete list for the forms of $\CC_k$, which includes special fibers with only multiplicity-1 components, see \cite{parametrizationofquadrics}.\\
\begin{center}
\begin{tabular}{|c|c|}
\hline
{$\CC_k$}&{Defining equations}\\
\hline
conic + double line &$x_1x_3=x_2^2+x_1x_4=0$\\
\hline
double conic& $x_1^2=x_2^2+x_3x_4=0$\\
\hline
double line + two lines&$ x_1^2+x_2^2=x_1x_3+\mu x_2x_4=0,\;\mu\in \{0,1\}$\\
\hline
triple line + line&$x_1x_2=x_1^2+x_2x_4=0$\\
\hline
two double lines&$x_2^2=x_1x_3+\mu x_2x_4=0,\;\mu\in \{0,1\}$\\
\hline
quadruple line&$x_1^2=x_2^2+ x_1x_3=0$\\
\hline
\end{tabular}
\end{center}

 Let the quadrics $F$ and $G$ be given by the following two polynomials respectively:
{\setlength\arraycolsep{2pt}
\begin{eqnarray}
\label{eq8}
 a_1x_1^2&+&a_2x_1x_2+a_3x_1x_3+a_4x_1x_4+a_5x_2^2+a_6x_2x_3+a_7x_2x_4+a_8x_3^2+a_9x_3x_4+a_{10}x_4^2,\nonumber\\
 b_1x_1^2&+&b_2x_1x_2+b_3x_1x_3+b_4x_1x_4+b_5x_2^2+b_6x_2x_3+b_7x_2x_4+b_8x_3^2+b_9x_3x_4+b_{10}x_4^2.\nonumber\\
 \end{eqnarray}}
 \begin{Proposition}
 \label{prop:normal4}
 Let $\CC$ be the $S$-scheme defined by the integral equation $\phi:F=G=0$, where $F$ and $G$ are given in (\ref{eq8}). Assume that $\CC_K$ is smooth.
 \begin{myitemize}
\item[(i)] If $\CC_k$ contains a multiplicity-$1$ component $\Gamma,$ then $\CC$ is normal at $\Gamma.$
\item[(ii)] If $\tilde{F}=x_1x_3$ and $\tilde{G}=x_2^2+x_1x_4,$ then $\CC$ is normal if and only if $$\nu(x_4F(0,0,x_3,x_4)-x_3G(0,0,x_3,x_4))=1.$$
\item[(iii)] If $\tilde{F}=x_1^2$ and $\tilde{G}=x_2^2+x_3x_4.$ Then $\CC$ is normal unless $F(0,x_2,x_3,x_4)\equiv\mu(x_2^2+x_3x_4)\mod t^2$, for some $\mu\in \OK$.
\item[(iv)] Assume that $\CC_k$ contains a line $\Gamma:\{x_1=x_2=0\}$ of multiplicity-$m,\;m\ge2,$ with $\tilde{F}=q(x_1,x_2)$ and $\tilde{G}=x_1x_3+\mu x_2x_4+q'(x_1,x_2)$, $\mu\in\{0,1\}$. If $\nu(F(0,0,x_3,x_4))=1,$ then $\CC$ is normal at $\Gamma$.
    \end{myitemize}
   \end{Proposition}
\begin{Proof}
$(i)$ Since $\CC_k$ is reduced at the generic point $\xi$ of $\Gamma$, we see that $\CC$ is normal at $\xi$.

Now we use Lemma \ref{lem:regularity} $(ii)$ to study the normality of $\CC$ at components of multiplicity greater than 1.
 The model $\CC$ is normal if and only if $F,G\not\in\mm_{\xi}^2$, and $F,G$ are linearly independent mod $\mm_{\xi}^2,$ for every generic point $\xi$ of $\CC_k.$ The linear independence condition is: For $\lambda_1,\lambda_2\in\OK[x_1,\ldots,x_4]_{\xi},$ if $\lambda_1F+\lambda_2G\in\mm_{\xi}^2,$ then $\lambda_1,\lambda_2\in\mm_{\xi}.$

$(ii)$ Let $\xi$ be the generic point of the double line $\{x_1=x_2=0\}$ in $\CC_k,$ then $\mm_{\xi}=\langle x_1,x_2,t\rangle.$ It is clear that $F,G\not\in\mm_{\xi}^2.$ If $\lambda_1F+\lambda_2G\in\mm_{\xi}^2,$ then the fact that $x_1$ and $t$ are linearly independent mod $\mm_{\xi}^2$ implies that $\lambda_1x_3+\lambda_2x_4\in\mm_{\xi},$ i.e., $\lambda_1\equiv \mu x_4 \mod\mm_{\xi}$ and $\lambda_2\equiv -\mu x_3\mod \mm_{\xi}$ for some $\mu\in\mathcal{O}_K.$ Thus $\CC$ is normal if and only if $\nu(f)=1,$ where $$f=x_4(a_8x_3^2+a_9x_3x_4+a_{10}x_4^2)-x_3(b_8x_3^2+b_9x_3x_4+b_{10}x_4^2).$$

$(iii)$ Let $\mm_{\xi}=\langle x_1, x_2^2+x_3x_4,t\rangle$ be the maximal ideal corresponding to the generic point $\xi$ of the conic. We have $G\not\in\mm_{\xi}^2.$ If $a_5x_2^2+a_9x_3x_4=tu(x_2^2+x_3x_4),u\in\mathcal{O}_K,$ then $F\not\in\mm_{\xi}^2$ if and only if $\nu(a_6x_2x_3+a_7x_2x_4+a_8x_3^2+a_{10}x_4^2)=1,$ otherwise $F\not\in\mm_{\xi}^2$ if and only if $\nu(a_5x_2^2+a_6x_2x_3+a_7x_2x_4+a_8x_3^2+a_9x_3x_4+a_{10}x_4^2)=1.$
If $\lambda_1F+\lambda_2G\in\mm_{\xi}^2,$ then $\lambda_2\in\mm_{\xi}.$ The reason is $t$ and $x_2^2+x_3x_4$ are linearly independent mod $\mm_{\xi}^2.$ Then the condition we obtained from $F\not\in\mm_{\xi}^2$ implies that $\lambda_1\in\mm_{\xi}.$

$(iv)$ Assume that $\xi$ is the generic point of $\Gamma:\{x_1=x_2=0\}$. The ideal $\mm_{\xi}$ is given by $\langle x_1,x_2,t\rangle.$ Since $\tilde{F}=q(x_1,x_2)$ and $\nu(a_8x_3^2+a_9x_3x_4+a_{10}x_4^2)=1,$ we have $F\not\in\mm_{\xi}^2.$
Since $\tilde{G}=x_1x_3+\mu x_2x_4+q'(x_1,x_2)$, we have $G\not\in\mm_{\xi}^2$ because $x_1x_3\not\in\mm_{\xi}^2$. Assume that $\lambda_1F+\lambda_2G\in\mm_{\xi}^2.$ Since $x_1,x_2$ and $t$ are linearly independent mod $\mm_{\xi}^2,$ it follows that $\lambda_2\in\mm_{\xi}.$ Moreover, as $\nu(a_8x_3^2+a_9x_3x_4+a_{10}x_4^2)=1,$ we get $\lambda_1\in\mm_{\xi}$.
  \end{Proof}

 \begin{Theorem}
\label{thm:minimal+soluble imlies normal}
Let $\phi$ be an integral genus one equation of degree $n,\;n\le 4,$ defining an $S$-scheme $\CC$. Assume that $\CC_K$ is smooth and $\CC_K(K)\ne\emptyset$.
If $\phi$ is minimal, then $\CC$ is normal.
\end{Theorem}
\begin{Proof}
If $\CC_k$ consists only of multiplicity-$1$ components, then $\CC_k$ is reduced and hence $\CC$ is normal. So, we only need to assume that $\phi$ is of degree $n,\;n\ge2,$ and $\CC_k$ contains a component of multiplicity greater than $1.$ Furthermore, we assume that $\CC$ is not normal, and hence $\phi$ does not satisfy the conditions included in Propositions \ref{prop:normal1,2,3} and \ref{prop:normal4}. Then we will prove that $\phi$ is not minimal by finding a genus one equation $K$-equivalent to $\phi$ whose discriminant has valuation less than $\Delta_{\phi}$. Recall that the discriminant varies by the 12th power of the determinant of the $K$-equivalence transformation, see \S \ref{sec:intro}.

 Let $n=2$ and $\phi:y^2+g(x)y=f(x)$. Since $\CC_k$ is a double line, we can assume that $\tilde{g}=\tilde{f}=0$. Since $\phi$ is not normal, we have $\nu(f)\ge 2.$ The equation $\phi$ is not minimal because it is $K$-equivalent to the equation $y^2+\frac{1}{t}g(x)y=\frac{1}{t^2}f(x).$

 Let $n=3$ and $\phi:F(x,y,z)=0$ as in Equation (\ref{eq7}). Now $\CC_k$ contains a multiplicity-$m$ component $\Gamma:\{y=0\}$, $m\ge2$. Since $\CC$ is not normal, we have $\nu(f_2),\nu(f_3)\ge1$ and $\nu(f_3)\ge 2$. Now $\phi$ is not minimal because it is $K$-equivalent to the equation $\frac{1}{t^2}F(x,ty,z)=0.$

 Let $n=4$ and $\phi:F=G=0,$ where $F$ and $G$ are given as in (\ref{eq8}). We will go through the different cases of Proposition \ref{prop:normal4}.

 Assume that $\CC_k:\{x_1x_3=x_2^2+x_1x_4=0\}$ and $\nu(x_4F(0,0,x_3,x_4)-x_3G(0,0,x_3,x_4))\ge2$. We use a matrix in $\GL_4(\OK)$ to get rid of the $x_1^2,x_1x_2$ and $x_1x_4$-terms in $F$ and of the $x_1^2,x_1x_2$ and $x_1x_3$-terms in $G$. We notice that in the equation
$$x_4F(0,0,x_3,x_4)-x_3G(0,0,x_3,x_4)=-b_8x_3^3+(a_8-b_9)x_3^2x_4+(a_9-b_{10})x_3x_4^2+a_{10}x_4^3,$$
we have $\min\{\nu(b_{8}),\nu(a_8-b_9),\nu(a_9-b_{10}),\nu(a_{10})\}\ge 2$. We apply the transformation $x_1\mapsto x_1-a_8x_3-a_9x_4,x_i\mapsto x_i,i=2,3,4,$ to get rid of the terms $a_8x_3^2$ and $a_9x_3x_4$. Thereafter, we obtain the genus one equation $\phi':F'=G'=0$, where
 {\setlength\arraycolsep{2pt}
 \begin{eqnarray}
 F'&=&x_1x_3+a_5x_2^2+a_6x_2x_3+a_7x_2x_4+a_{10}x_4^2,\nonumber\\
 G'&=&x_1x_4+x_2^2+b_6x_2x_3+b_7x_2x_4+b_8x_3^2+(b_9-a_8)x_3x_4+(b_{10}-a_9)x_4^2.\nonumber
 \end{eqnarray}}
 We deduce that $\phi'$ is not minimal because it is $K$-equivalent to the equation $$\frac{1}{t^2}F'(t^2x_1,tx_2,x_3,x_4)=\frac{1}{t^2}G'(t^2x_1,tx_2,x_3,x_4)=0.$$

 Assume that $\CC_k:\{x_1^2=x_2^2+x_3x_4=0\},$ and that $\nu(F(0,x_2,x_3,x_4)-\mu(x_2^2+x_3x_4))\ge2$ for some $\mu\in\OK$. Then $\phi$ is not minimal because it is $K$-equivalent to
 $$\frac{1}{t^2}(F(tx_1,x_2,x_3,x_4)-\mu  G(tx_1,x_2,x_3,x_4))=G(tx_1,x_2,x_3,x_4)=0.$$

Now assume that $\CC_k$ contains a line $\Gamma:\{x_1=x_2=0\}$ of multiplicity-$m,\;m\ge2,$ with $\tilde{F}=q(x_1,x_2)$ and $\tilde{G}=x_1x_3+\mu x_2x_4+q'(x_1,x_2),$ where $\mu\in\{0,1\}.$ Assume that $\nu(a_8x_3^2+a_9x_3x_4+a_{10}x_4^2)\ge2$. Then $\phi$ is not minimal because it is $K$-equivalent to
$$\frac{1}{t^2}F(tx_1,tx_2,x_3,x_4)=\frac{1}{t}G(tx_1,tx_2,x_3,x_4)=0.$$
\end{Proof}

\section{Criteria for minimality}
\label{sec:criteria}

We state the first main theorem of this paper.

\begin{Theorem}
\label{thm:minimality}
Let $\phi$ be an integral genus one equation of degree $n,\;n\le4$. Assume moreover that $\phi$ defines a normal $S$-scheme $\CC$, and that $\CC_K$ is smooth with $\CC_K(K)\ne\emptyset.$
Then $\phi$ is minimal, see Definition \ref{def:integral-minimal}, if and only if $\phi$ is geometrically minimal, see Definition \ref{def:geometrically minimal}.
\end{Theorem}

Theorem \ref{thm:minimality} is known for the case $n=1,$ see (\cite{Liubook}, \S 9.4) or \cite{Brian}. When $n=2$, Liu proved that if $\phi$ is minimal, then $\phi$ is geometrically minimal, see (\cite{LiuModeles}, Corollaire 5 (a)). We introduce a proof which works for $n$ when $n\le4.$

\subsection{Canonical sheaves of $S$-models}

Let $E$ be an elliptic curve with minimal proper regular model $E^{min}$. If $\CC$ is an $S$-model for $E,$ then the canonical sheaf $\omega_{\CC/S}$ of $\CC$ satisfies $\omega_{\CC/S}|_E=\omega_{E/K}$, moreover the restriction of the canonical sheaf $\omega_{\CC/S}$ on $E$ gives a canonical injection $H^0(\CC,\omega_{\CC/S})\hookrightarrow H^0(E,\omega_{E/K})$. We have $H^0(\CC,\mathcal{O}_{\CC})=\OK.$ In particular, if $\omega_{\CC/S}=\omega\mathcal{O}_{\CC},$ then $H^0(\CC,\omega_{\CC/S})=\omega\OK.$ In addition, there exists an $\omega_0\in H^0(E^{min},\omega_{E^{min}/S})$ such that $\omega_{E^{min}/S}= \omega_0\mathcal{O}_{E^{min}},$ see (\cite{Brian}, Example 7.7).

\begin{Lemma}
\label{lem:H0Emin H0Cdelta}
Let $E$ be an elliptic curve over $K$ with minimal proper regular model $E^{min}$. Let $\CC$ be a normal $S$-scheme with $\CC_K\isom E$. Let $\Cdes\to\CC$ be the minimal desingularisation of $\CC$. Then the following statements hold.
\begin{myitemize}
\item[(i)] $H^0(E^{min},\omega_{E^{min}/S})=H^0(\Cdes,\omega_{\Cdes/S})\subseteq H^0(\CC,\omega_{\CC/S})$
as subgroups in $H^0(E,\omega_{E/K}).$
\item[(ii)] If $\Cdes\isom E^{min}$, then $H^0(E^{min},\omega_{E^{min}/S})=H^0(\CC,\omega_{\CC/S}).$
\end{myitemize}
\end{Lemma}
\begin{Proof}
$(i)$ The equality holds because $\widetilde{\CC}$ and $E^{min}$ are two regular $S$-models for $E$, see (\cite{Liubook}, Corollary 9.2.25 (b)). The inequality holds because $\CC$ is obtained from $\Cdes$ as a contraction of a finite number of irreducible components, see (\cite{Liubook}, Lemma 9.2.17 (a)).

$(ii)$ Since $\Cdes\isom E^{min}$, we have a contraction morphism $f:E^{min}\to \CC$. Therefore, $f_*\omega_{E^{min}/S}=\omega_{\CC/S}$, see (\cite{Liubook}, Corollary 9.4.18 (b)).
\end{Proof}

\begin{Proposition}
\label{prop:canoniacal sheaf of degree n model}
Let $\CC$ be an $S$-scheme defined by an integral genus one equation $\phi$ of degree $n$. Assume that $\CC_K$ is smooth. Then $\omega_{\CC/S}=\omega_{\phi} \mathcal{O}_{\CC},$ where $\omega_{\phi}\in H^0(\CC_K,\omega_{\CC_K/K})$ is
\begin{myitemize}
\item[(i)] if $n=1$: $\omega_{\phi}= du/(2v+a_1u+a_3),\textrm{ where }u=x/z,v=y/z\in K(\CC),$
\item[(ii)] if $n=2$: $\omega_{\phi}= dx/(2y+g(x)),$
\item[(iii)] if $n=3$: $\omega_{\phi}=du/(\partial{F}/\partial{v}),\textrm{ where } u=x/z, v=y/z\in K(\CC), $
\item[(iv)] if $n=4$: $\omega_{\phi}= d u_2/  (\frac{\partial{F_1}}{\partial{u_4}}\frac{\partial{F_2}}{\partial{u_3}}-\frac{\partial{F_1}}{\partial{u_3}}\frac{\partial{F_2}}{\partial{u_4}}), \textrm{ where }u_i=x_i/x_1 \in K(\CC),\;i=2,3,4.$
\end{myitemize}
\end{Proposition}
\begin{Proof}
This is a direct consequence of (\cite{Liubook}, Corollary 6.4.14).
\end{Proof}

\begin{Lemma}
\label{lem:rel between delta.omega}
Let $\phi_1,\phi_2$ be two $K$-equivalent integral genus one equations of degree $n,\;n\le4$. Assume that $\phi_i$ defines an $S$-scheme $\CC_i$ with $C_i:=(\CC_i)_K$ being smooth. If $\omega_{\CC_i/S}=\omega_{\phi_i}\mathcal{O}_{\CC_i}$, then
$$\Delta_{\phi_1}\omega_{\phi_1}^{\otimes 12}=\lambda\Delta_{\phi_2} \omega_{\phi_2}^{\otimes 12}\in H^0(C_1,\omega_{C_1/K})^{\otimes 12}=H^0(C_2,\omega_{C_2/K})^{\otimes12},\textrm{ where }\lambda\in\mathcal{O}^*_K.$$
\end{Lemma}
\begin{Proof}
Assume that $\phi_1=g.\phi_2$ where $g$ is a transformation defining the $K$-equivalence.
 The transformation $g$ defines an isomorphism $\gamma:C_1\isom C_2$ which satisfies $\gamma^*\omega_{\phi_2}=(\det g) \omega_{\phi_1},$ see (\cite{FiInvariants}, Proposition 5.19). Hence, $\omega_{\phi_2}=\alpha(\det g) \omega_{\phi_1}$ as elements in $H^0(C_i,\omega_{C_i/K}),$ where $\alpha\in\mathcal{O}^*_K.$ Recall that $\Delta_{\phi_1}=(\det g)^{12}. \Delta_{\phi_2}.$ It follows that $\Delta_{\phi_1}\omega_{\phi_1}^{\otimes 12}=\alpha^{12}\Delta_{\phi_2} \omega_{\phi_2}^{\otimes 12}.$
\end{Proof}

If $\phi_1$ is minimal, then we call the integer $m$ such that $\omega_{\phi_2}=ut^{-m}\omega_{\phi_1},u\in\mathcal{O}^*_K,$ the {\em level} of $\phi_2,$ and denote it by $\level(\phi_2).$ The above corollary implies that the level of an integral genus one equation of degree $n$ does not depend on the choice of the minimal genus one equation $\phi_1.$
Notice that $\nu(\Delta_{\phi_2})=\nu(\Delta_{\phi_1})+12 \level(\phi_2).$

\begin{Lemma}
\label{lem:H0 subset H0}
Keep the hypothesis of Lemma \ref{lem:rel between delta.omega}. Then we have
$H^0(\CC_1,\omega_{\CC_1/S})\subseteq H^0(\CC_2,\omega_{\CC_2/S})$ as sub-$\OK$-modules of $H^0(C_i,\omega_{C_i/K})$ if and only if $\nu(\Delta_{\phi_1})\le\nu(\Delta_{\phi_2}).$ Moreover, the equality of the two submodules holds if and only if $\phi_1$ and $\phi_2$ have the same level.
\end{Lemma}
\begin{Proof}
 The assumption $H^0(\CC_1,\omega_{\CC_1/S})\subseteq H^0(\CC_2,\omega_{\CC_2/S})$ is equivalent to $\omega_{\phi_1}\in\omega_{\phi_2}\OK.$ Thus Lemma \ref{lem:rel between delta.omega} asserts that $\Delta_{\phi_2}\in\Delta_{\phi_1}\OK,$ i.e., $\nu(\Delta_{\phi_1}))\le\nu(\Delta_{\phi_2}).$

The equality of the sub-$\OK$-modules $H^0(\CC_1,\omega_{\CC_1/S})= H^0(\CC_2,\omega_{\CC_2/S})$ means that $\omega_{\phi_1}\OK=\omega_{\phi_2}\OK$ as $\OK$-modules, i.e., $\omega_{\phi_1}\in\omega_{\phi_2}\mathcal{O}^*_K$. Hence, $\phi_1$ and $\phi_2$ have the same level by Lemma \ref{lem:rel between delta.omega}.
\end{Proof}

\subsection{Constructing genus one equations}
\label{sec:constructing minimal degree n models}

Let $E$ be an elliptic curve over $K$ with
minimal proper regular model $E^{min}$. Let $P\in E(K)$. We will denote the Zariski closure of $\{P\}$ in $E^{min}$ by $\overline{\{P\}}$.

When $n=1,$ set $D_1=3.\overline{\{P\}}$ where $P\in E(K).$
When $n\ge2,$ let $\sum(P_i)\in\Div(E)$ be a $K$-rational divisor on $E$ of degree $n.$ Assume moreover that $\overline{\{P_i\}}\cap E^{min}_k$ is contained in one and only one irreducible component of $E^{min}_k.$ Consider the divisor $D_n$ on $E^{min}$, and the $S$-model $\CC_n$ for $E$ given by
$$D_n=\sum\overline{\{P_i\}},\textrm{ and }\CC_n:=\Proj(\bigoplus_{m=0}^{\infty}H^0(E^{min},\mathcal{O}_{E^{min}}(mD_n))).$$
 There is a canonical morphism $u_n:E^{min}\rightarrow \CC_n$ contracting all the irreducible components of $E^{min}_k$ except the ones having nonempty intersection with $D_n$, see Theorem \ref{th:contraction}.

\begin{Lemma}
\label{lem:free module}
Let $D_n, \;n\in\{1,2,3,4\},$ be as above. Then $H^0(E^{min},\mathcal{O}_{E^{min}}(mD_n)),m\ge1,$ is a free $\OK$-module of rank $3m$ if $n=1,$ and of rank $mn$ if $n\ge2.$
\end{Lemma}
\begin{Proof}
It is known that $H^0(E^{min},\mathcal{O}_{E^{min}}(mD_n))\otimes_{\OK}K\isom H^0(C,\mathcal{O}_C(mD_n|_C))$, see for example (\cite{Liubook}, Corollary 5.2.27). By virtue of the Riemann-Roch Theorem, $H^0(C,\mathcal{O}_C(mD_n|_C))$ is a $3m$-dimensional $K$-vector space when $n=1,$ and an $mn$-dimensional $K$-vector space when $n\ge2.$

Since $\mathcal{O}_{E^{min}}(mD_n)$ is an invertible sheaf on $E^{min}$, it follows that $H^0(E^{min},\mathcal{O}_{E^{min}}(mD_n))$ is a flat $\OK$-module. Therefore,
$H^0(E^{min},\mathcal{O}_{E^{min}}(mD_n))$ is a finitely generated flat module over the local ring $\OK$, hence it is free.
\end{Proof}

\begin{Theorem}
\label{thm:C_n is minimal}
Let $\CC_n$ and $D_n,\;n\le4,$ be as above. Then there exists an integral genus one equation $\phi_n$ of degree $n$ defining $\CC_n$.
\end{Theorem}
\begin{Proof}
  For $ n=1,$ see (\cite{Liubook}, \S 9.4).
For $n\ge2$, we pick a basis $\{x_1,\ldots,x_n\}$ of $H^0(E^{min},\mathcal{O}_{E^{min}}(D_n))$. Consider the morphism
$\lambda_n:E^{min}\longrightarrow \PP_S^{n-1}$
associated to the basis $\{x_1,\ldots,x_n\}$.

 For $n=2$, put $x_1=x$ and $x_2=1$. We have $\PP_S^1=\Spec\OK[x]\cup\Spec\OK[1/x]$. Let $U=\lambda_2^{-1}(\Spec\OK[x])$ and $V=\lambda_2^{-1}(\Spec\OK[1/x])$. We have $\CC_2=U\cup V.$ Taking the integral closure of $\OK[x]$ in $K(\CC_2)$, we have $$\mathcal{O}_{\CC_2}(U)=\OK[x]\oplus y\OK[x], \textrm{ for some }y\in\mathcal{O}_{\CC_2}(U),$$
moreover there exist $g(x),f(x)\in\OK[x]$ such that $\deg g\le 2,\;\deg f \le 4$ and $y^2+g(x)y=f(x)$, see (\cite{LiuModeles}, Lemme 1). Following the same argument for $\mathcal{O}_{\CC_2}(V)$, we deduce that $\CC_2$ is the union of the two affine open schemes
\begin{eqnarray}
U=\Spec\OK[x,y]/(y^2+g(x)y-f(x)),\textrm{ }V=\Spec\OK[w,z]/(z^2+w^2 g(1/w)z-w^4f(1/w)),\nonumber
\end{eqnarray}
where $w=1/x,\;z=y/x^2.$ Hence, we are done when $n=2$.

Now let $n=3,4$. Let $Z_n$ be the closed subset $\lambda_n(E^{min})\subset\PP^{n-1}_S$ endowed with the reduced scheme structure. According to the description of the contraction morphism included in the proof of (\cite{Liubook}, Proposition 8.3.30), the morphism $\lambda_n:E^{min}\to Z_n\subseteq \PP_S^{n-1}$ factors into $u_n:E^{min}\to \CC_n$ followed by $v_n:\CC_n\to Z_n,$
where $v_n$ is the normalisation morphism. It is understood that $v_n$ is a finite morphism, hence for an irreducible component $\Gamma$ of $E^{min}_k$,
 $\lambda_n(\Gamma)$ is a point if and only if $u_n(\Gamma)$ is a point.
 In other words, the special fibers of $\CC_n$ and $Z_n$ have the same number of irreducible components. We are going to show that $Z_n$ is defined by an integral genus one equation of degree $n$. Then we show that $\CC_n\isom Z_n$.

When $n=3$, the free $\OK$-module $H^0(E^{min},\mathcal{O}_{E^{min}}(3D_3))$ is of rank-$9$, see Lemma \ref{lem:free module}, but it contains the $10$ elements $x_1^3, x_2^3,x_3^2, x_1^2x_2, x_1^2x_3, x_2^2x_1, x_2^2x_3, x_3^2x_1, x_3^2x_2, x_1x_2x_3$. It follows that there are $a_i\in\OK$ such that
$$F:=a_1x_1^3+a_2x_2^3+a_3x_3^3+a_4x_1^2x_2+a_5x_1^2x_3+a_6x_2^2x_1+a_7 x_2^2x_3+a_8 x_3^2x_1+a_9 x_3^2x_2+a_{10}x_1x_2x_3=0.$$
Rescaling $x_1,\;x_2$ and $x_3$, we can assume that there is at least one $a_i\in\mathcal{O}^*_K$. Now $Z_3$ is contained in $\Proj\OK[x_1,x_2,x_3]/(F)$.

When $n=4$, we consider the $10$ elements $x_1^2,x_1x_2,x_1x_3,x_1x_4,x_2^2,x_2x_3,x_2x_4,x_3^2,x_3x_4,x_4^2$ in the rank-$8$ free $\OK$-module $H^0(E^{min},\mathcal{O}_{E^{min}}(2D_4))$. They satisfy two quadrics $Q$ and $R$ with coefficients in $\OK$. Moreover, by rescaling $x_1,x_2,x_3$ and $x_4$ we can assume that both $Q$ and $R$ have at least one coefficient in $\mathcal{O}_K^*$. Now $Z_4$ is contained in the intersection of $Q$ and $R$.

We want to show that $Z_n= \Proj\OK[x_1,\ldots,x_n]/I_n$, where $I_3=(F)$ and $I_4=(Q,R)$. Recall that both schemes are of dimension 2.
 Since $Z_n\subseteq\Proj\OK[x_1,\ldots,x_n]/I_n$, we have $\Proj\OK[x_1,\ldots,x_n]/I_n=Z_n\cup Z_n'$, for some closed subscheme $Z_n'\subset\PP^{n-1}_S,$ where $Z_n'\ne \Proj\OK[x_1,\ldots,x_n]/I_n$. Since  $\Proj(\OK[x_1,\ldots,x_n]/I_n\otimes K)$ is irreducible, $\Proj\OK[x_1,\ldots,x_n]/I_n$ is irreducible itself. It follows from the definition of irreducibility that $Z_n'=\emptyset$, and the closed subscheme $Z_n$ is $\Proj\OK[x_1,\ldots,x_n]/I_n$.

 Now both $\CC_n$ and $Z_n$, $n=3,4,$ have dimension $2$,
 their generic fibers are isomorphic, and their special fibers  have the same number of irreducible components.
 By virtue of (\cite{Liubook}, Remark 8.3.25), $v_n:\CC_n\to Z_n$ is an isomorphism.
\end{Proof}

\begin{Remark}
\label{rem1}
 Let $n\ge2.$ Let $\phi$ be an integral genus one equation of degree $n$ defined by a hyperplane section $H.$ Assume moreover that $\phi$ defines a morphism $E\to\PP^{n-1}_K$. Then we can choose $D_n=(n-1)\overline{\{P\}}+\overline{\{Q\}}$ where $P,Q\in E(K)$ are such that $(n-1)P+Q\sim H.$ It follows that $\phi_n$ is $K$-equivalent to $\phi$.
\end{Remark}


\subsection{Proof of Theorem \ref{thm:minimality}}

\begin{ProofOf}{Theorem \ref{thm:minimality}}
We first assume that $\phi$ is geometrically minimal, i.e., that $\widetilde{\CC}\isom E^{min}$. Thus we have a contraction morphism $f : E^{min}\rightarrow \CC.$
We claim that if $\CC'$ is a normal $S$-scheme defined by an integral genus one equation $\phi'$ which is $K$-equivalent to $\phi$, then
$H^0(\CC,\omega_{\CC/S})\subseteq H^0(\CC',\omega_{\CC'/S})$,
and hence $\nu(\Delta_{\phi})\le\nu(\Delta_{\phi'}),$ see Lemma \ref{lem:H0 subset H0}. Therefore, $\phi$ is minimal.

Lemma \ref{lem:H0Emin H0Cdelta} $(i)$ shows that $H^0(E^{min},\omega_{E^{min}/S})\subseteq  H^0(\CC',\omega_{\CC'/S}).$
The fact that $\CC$ is obtained from $E^{min}$ by contraction implies that
$H^0(\CC,\omega_{\CC/S}) = H^0(E^{min},\omega_{E^{min}/S})$, see Lemma \ref{lem:H0Emin H0Cdelta} $(ii)$. Thus the claim is proved.

Before we proceed with the proof of the second part of the theorem we need the following lemma.

\begin{Lemma}
\label{lem:H0 Emin and H0 C when min}
Let $\phi$ be a genus one equation of degree $n$ defining a normal $S$-scheme $\CC$ with $\CC_K$ being smooth and $\CC_K(K)\ne\emptyset.$ Let $E^{min}$ be the minimal proper regular model of the Jacobian of $\CC_K$. Then $H^0(E^{min},\omega_{E^{min}/S})= H^0(\CC,\omega_{\CC/S})$ as sub-$\OK$-modules of $H^0(\CC_K,\omega_{\CC_K/K})$ if and only if $\phi$ is minimal.
\end{Lemma}
\begin{Proof}
Assume that $H^0(E^{min},\omega_{E^{min}/S})=H^0(\CC,\omega_{\CC/S}).$ Let $\CC'$ be a normal $S$-scheme defined by a genus one equation $K$-equivalent to $\phi$. By virtue of Lemma \ref{lem:H0Emin H0Cdelta} $(i)$, we have
$H^0(\CC,\omega_{\CC/S})\subseteq H^0(\CC',\omega_{\CC'/S})$.
Therefore, $\phi$ is minimal by Lemma \ref{lem:H0 subset H0}.

Now assume that $\phi$ is minimal. According to Theorem \ref{thm:C_n is minimal} and Remark \ref{rem1}, there exists an integral genus one equation $\phi'$ which is $K$-equivalent to $\phi$. Moreover, $\phi'$ is geometrically minimal because the minimal desingularisation of the $S$-scheme $\CC'$ defined by $\phi'$ is isomorphic to
$E^{min}.$ Therefore, $H^0(E^{min},\omega_{E^{min}/S})=H^0(\CC',\omega_{\CC'/S})$ by Lemma \ref{lem:H0Emin H0Cdelta} $(ii)$. Moreover, $\phi'$ is minimal by the first part of Theorem \ref{thm:minimality}.
Since the genus one equations $\phi$ and $\phi'$ are both minimal, in particular they have the same level, Lemma \ref{lem:H0 subset H0} implies that
 $H^0(\CC',\omega_{\CC'/S})=H^0(\CC,\omega_{\CC/S})$, and we are done.
 \end{Proof}

Now we conclude the proof of Theorem \ref{thm:minimality}.

Assume that $\phi$ is minimal.
We assume on the contrary that $\Cdes\not\isom E^{min},$ and therefore $\CC_k$ contains an exceptional divisor $\Gamma.$
By (\cite{Liubook}, Proposition 9.3.10), we have $\deg \omega_{\Cdes/S}|_{\Gamma}<0$. It follows that $H^0(\Gamma,\omega_{\Cdes/S}|_{\Gamma})=0,$ therefore $\omega_{\Cdes/S}$ is not generated by its global sections on $\Gamma$. But we have
    $$\omega_{\Cdes/S}|_{\Gamma}=\omega_{\CC/S}|_{\Gamma}=\omega_{\phi}\mathcal{O}_{\CC}|_{\Gamma},$$
     where $\omega_{\phi}$ is given as in Proposition \ref{prop:canoniacal sheaf of degree n model}. The global sections of $\omega_{\Cdes/S}$ are $$H^0(\Cdes,\omega_{\Cdes/S})=H^0(E^{min},\omega_{E^{min}/S})=H^0(\CC,\omega_{\CC/S})=\omega_{\phi}\OK,$$ where the second equality is justified by $\CC$ being minimal, see
 Lemma \ref{lem:H0 Emin and H0 C when min}. Therefore,
    $\omega_{\Cdes/S}$ is generated by its global sections at every $x\in\Gamma,$ which is a contradiction.
\end{ProofOf}

\begin{Corollary}
\label{cor:critieria for minimality}
Let $\phi,\CC$ and $\Cdes$ be as in Theorem \ref{thm:minimality}, and $\omega_{\phi}$ as in Proposition \ref{prop:canoniacal sheaf of degree n model}. Then the following are equivalent.\\
\begin{tabular}{p{7cm}p{8cm}}
(i) $\phi$ is minimal & (ii) $\phi$ is geometrically minimal \\
(iii) $\omega_{\Cdes/S}=\omega_{\phi}\mathcal{O}_{\Cdes}$ & (iv) $H^0(\Cdes,\omega_{\Cdes/S})=H^0(\CC,\omega_{\CC/S})=\omega_{\phi}\OK.$
\end{tabular}
\end{Corollary}

\section{Existence of global models}
\label{sec:global models}

In order to prove the global result included in Theorem \ref{thm1}, we need the following local result.


\begin{Lemma}
\label{lem:ord delta}
Let $\phi$ be a minimal genus one equation of degree $n$ defining a normal $S$-scheme $\CC$ such that $\CC_K$ is smooth and $\CC_K(K)\ne\emptyset.$
Let $E$ be the Jacobian of $\CC_K$ with minimal discriminant $\Delta$. Then $\nu(\Delta_{\phi})=\nu(\Delta)$.
\end{Lemma}
\begin{Proof}
Let $\EE$ and $E^{min}$ be the minimal Weierstrass model and minimal proper regular model of $E$ respectively.
We identify $\CC_K$ and $E$ via an isomorphism $\gamma:\CC_K\isom_K E$. For explicit formulae for $\gamma$, see (\cite{FiInvariants}, \S 6) or (\cite{FiStCr}, \S 2). Let $\omega_{\phi}$ and $\omega$ be the generators of the canonical sheaves $\omega_{\CC/S}$ and $\omega_{\EE/S}$ given as in proposition \ref{prop:canoniacal sheaf of degree n model}.

Now according to Lemma \ref{lem:H0 Emin and H0 C when min}, we have $H^0(E^{min},\omega_{E^{min}/S})=H^0(\EE,\omega_{\EE/S})=\omega\OK$, and $H^0(E^{min},\omega_{E^{min}/S})=H^0(\CC,\omega_{\CC/S})=\omega_{\phi}\OK$ as sub-$\OK$-modules of $H^0(E,\omega_{E/K})$. Therefore, $\omega_{\phi}=\lambda\omega$ for some $\lambda\in\mathcal{O}_K^*.$ According to Lemma \ref{lem:rel between delta.omega}, we have $\Delta_{\phi}.\omega_{\phi}^{\otimes12}=\Delta.\omega^{\otimes12}$ up to a unit. Hence, $\Delta_{\phi}=\Delta$ up to a unit in $\mathcal{O}_K^*$. In other words, $\nu(\Delta_{\phi})=\nu(\Delta).$
\end{Proof}

\begin{Corollary}
\label{cor:delta and delta min}
Let $\phi$ be an integral genus one equation defining a normal $S$-scheme $\CC$ such that $\CC_K$ is smooth and $\CC_K(K)\ne\emptyset.$ Let $\Cdes\to\CC$ be the minimal desingularisation of $\CC$. Let $E$ be the Jacobian of $\CC_K$ with minimal discriminant $\Delta$. Then  $$\nu(\Delta_{\phi})=\nu(\Delta)+12\length_{\OK}(H^0(\CC,\omega_{\CC/S})/H^0(\Cdes,\omega_{\Cdes/S})),$$
where $\length_{\OK}$ denotes the length of an $\OK$-module.
\end{Corollary}
\begin{Proof}
  Theorem \ref{thm:C_n is minimal} and Remark \ref{rem1} admit the existence of an integral geometrically minimal genus one equation $\phi'$ $K$-equivalent to $\phi$. Moreover, Theorem \ref{thm:minimality} implies that $\phi'$ is minimal. We have $\nu(\Delta_{\phi'})=\nu(\Delta)$ by Lemma \ref{lem:ord delta}. We only need to show that $\nu(\Delta_{\phi})=\nu(\Delta_{\phi'})+12\length_{\OK}(H^0(\CC,\omega_{\CC/S})/H^0(\Cdes,\omega_{\Cdes/S}))$. The latter follows immediately from the fact that $\Delta_{\phi}.\omega_{\phi}^{\otimes 12}=\Delta_{\phi'}.\omega_{\phi'}^{\otimes 12}$ up to a unit, see Lemma \ref{lem:rel between delta.omega}, and $H^0(\CC',\omega_{\CC'/S})=H^0(\Cdes,\omega_{\Cdes/S})\subseteq H^0(\CC,\omega_{\CC/S})=\omega_{\phi}\mathcal{O}_{\CC}$, where $\CC'$ is the normal $S$-scheme defined by $\phi'$.
\end{Proof}

We observe that $\length_{\OK}(H^0(\CC,\omega_{\CC/S})/H^0(\Cdes,\omega_{\Cdes/S}))$ in the above corollary is an interpretation for $\level(\phi)$ mentioned after Lemma \ref{lem:rel between delta.omega}.
Now we conclude with the proof of Theorem \ref{thm1}.

\begin{ProofOf}{Theorem \ref{thm1}}
We write $\mathcal{O}_{\nu}$ for the ring of integers of the completion $F_{\nu}$. Set $S_{\nu}=\Spec\mathcal{O}_{\nu}$.

 By virtue of Theorem \ref{thm:C_n is minimal} and Remark \ref{rem1}, $\phi$ is $F_{\nu}$-equivalent to a geometrically minimal $\mathcal{O}_{\nu}$-integral genus one equation $\phi_{\nu}$. Moreover, $\phi_{\nu}$ is minimal, see Theorem \ref{thm:minimality}. According to Lemma \ref{lem:ord delta}, $\nu(\Delta_{\phi_{\nu}})=\nu(\Delta)$ for every non-archimedean place $\nu$.

 Since $F$ is of class number $1$, we can use strong approximation to find an $\mathcal{O}_F$-integral genus one equation $\phi'$ which is $F$-equivalent to $\phi$ and $\nu(\Delta_{\phi'})=\nu(\Delta_{\phi_{\nu}})=\nu(\Delta)$, for every non-archimedean place $\nu$. See \cite{Fiminimise} for details when $n=2,3$. The case $n=4$ is similar. Since the discriminant is of weight 12, we have $\Delta_{\phi'}=\lambda^{12}\Delta$ for some $\lambda\in\mathcal{O}_F^*.$ By rescaling the coefficients of the defining polynomials of $\phi'$, we can assume that $\lambda=1.$
\end{ProofOf}

\hskip-18pt\emph{\bf{Acknowledgements.}}
This paper is based on the author's Ph.D. thesis \cite{SadekThesis} at Cambridge University. The author would like to thank his supervisor Tom Fisher for guidance and many useful comments on the manuscript.

\bibliographystyle{plain}
\footnotesize
\bibliography{thesisreferences}

\end{document}